\magnification=1200
\input amstex
\documentstyle{amsppt}

\pageheight {49pc}
\vcorrection{-2pc}


\def\leftheadline {\line{\hfill\eightrm MARSHALL A. WHITTLESEY\hfill\folio}}
\def\rightheadline{\line{\hfill \eightrm POLYNOMIAL HULLS AND AN
OPTIMIZATION PROBLEM\hfill \folio}}

\headline{\ifodd \pageno\rightheadline \else\leftheadline\fi}
\def\C{\hbox{\bf C$\,$}}
\def\R{\hbox{\bf R$\,$}}
\def\opdisk{\hbox {\rm int$\,\Delta\,\,$}}

\def\origin{\buildrel\rightarrow\over 0}
\def\ui{\hbox{$u_1$}}
\def\uo{\hbox{$u_2$}}

\def\uk{\hbox{$u_3$}}
\def\ul{\hbox{$u_5$}}
\def\um{\hbox{$u_4$}}
\def\un{\hbox{$u_6$}}

\def\uq{\hbox{$u_8$}}
\def\us{\hbox{$u_9$}}
\def\ur{\rho}

\def\monge{\tau_z }

\def\hullthm{1\,}
\def\optthm{2\,}
\def\fathull{3\,}

\def\lemcenter{1\,}
\def\lemsup{2\,}
\def\lemdelta{3\,}
\def\lemmollify{4\,}
\def\lemlocmollify{5\,}
\def\lemdeffuncA{6\,}
\def\lemdeffuncB{7\,}

\topmatter
\title  Polynomial hulls and an optimization problem 
\endtitle
\author Marshall A. Whittlesey
\endauthor
\affil  California State University, San Marcos
\endaffil
\address  Department of Mathematics, California State University
San Marcos, San Marcos CA 92096
\endaddress
\email mwhittle\@csusm.edu
\endemail

\abstract We say that a subset of $\C^n $ is hypoconvex if its
complement is the union of complex hyperplanes.  We say it is strictly
hypoconvex if it is smoothly bounded hypoconvex and at every point of
the boundary the real Hessian of its defining function is positive
definite on the complex tangent space at that point.  Let $B_n$ be
the open unit ball in $\C^n$.  Suppose $K$ is
a $C^\infty$ compact manifold in $\partial B_1\times\C^n$, $n>1$, 
diffeomorphic to $\partial B_1\times \partial B_n$,
each of whose fibers $K_z$ over $\partial B_1$ bounds a strictly hypoconvex connected open set.  Let
$\widehat K$ be the polynomial hull of $K$.  Then we show that
$\widehat K\setminus K$ is the union of graphs of analytic vector
valued functions on $B_1$.  This result shows that an unnatural
assumption regarding the deformability of $K$ in an earlier version of
this result is unnecessary.  Next, we study an $H^\infty$ optimization
problem.  If $\rho $ is a $C^\infty$ real-valued function on
$\partial B_1\times\C^n$, we show that the infimum $\gamma_\rho=
\inf_{f\in H^\infty(B_1)^n}
\|\rho(z,f(z))\|_\infty$ is attained by a unique bounded $f$
provided that the set
$\{(z,w)\in\partial B_1\times\C^n|\rho(z,w)\leq\gamma_\rho\}$ has bounded connected strictly hypoconvex fibers over the circle.
\endabstract

\thanks 
\endthanks

\endtopmatter

\hsize=6.75truein
\hoffset=-.15in
\voffset=0.3 true in
\vsize=8.7truein

\document
\nopagenumbers
\NoBlackBoxes

$\S 1$ {\bf Introduction and results.}

The purpose of this work is to strengthen the results of [Wh] and [V]
regarding the existence of analytic structure in a class of polynomial
hulls and the solution of an $H^\infty $ optimization problem.

For $w\in\C^n$ let $w_j$ denote the $j^{th}$ coordinate of $w$ and let
$|w|$ denote the standard Euclidean norm.  If $Y$ is a compact set in $\C^n$, then
the {\it polynomial (convex) hull} $\widehat Y$ of $Y$ is given by
$$\hbox {$\widehat Y =
\{z\in\C^n\bigl ||P(z)|\leq
{\displaystyle \sup _{w\in Y}} |P(w)|$ for all polynomials $P$ on
\C$^n\}$.}$$ A basic question concerning polynomial hulls is the
following: given $Y$ and $y\in \widehat Y\setminus Y$, does there
exist an analytic variety through $y$ with boundary in $Y$?  If so,
that variety gives a reason why $y$ is in $\widehat Y$, by virtue
of the maximum modulus principle for analytic varieties.

Let $\Delta $ denote the closed unit disk in \C and $\Gamma $ its
boundary.  Let $B_n$ be the open unit ball in $\C^n$ and $\origin$ the origin in
$\C^n$.  Let $K\subset \Gamma\times\C^n$.  Does $\widehat K$ contain
analytic varieties?  Does $\widehat K$ contain graphs of analytic
functions over the unit disk?  Answers to these questions have been formulated, 
in works by several authors, in terms of the fibers $K_z\equiv\{w\in\C^n:(z,w)\in
K\}$.  In [AW,S2], it was shown that if the $K_z$ are convex then
$\widehat K
\setminus K$ is the union of analytic graphs over the disk.
The same was shown in the case $n=1$ for fibers $K_z$ which are
connected and simply connected in [F] and [S1].  See [\v Ce2] for a theorem 
in the case  where the fibers are closures of completely circled pseudoconvex 
domains.  Examples of $K$
where $\widehat K$ does not contain analytic graphs over the disk have
been constructed in [HMe] and [\v Ce1].  In [HMa], Helton and Marshall
conjecture that if the fibers $K_z$ are (connected and) possess
a property called {\it hypoconvexity} then the set $\widehat K\setminus
K$ is the union of analytic graphs over the open disk.

We say that a subset of $\C^n$ is 
{\it hypoconvex} if its complement is the union of complex hyperplanes
(of complex dimension $n-1$); {\it strictly hypoconvex} if it is
smoothly bounded hypoconvex and the real Hessian of its defining
function is positive definite on the complex tangent space at every
point in the boundary.  (The terminology in the literature on this
subject varies; instead of hypoconvex, the reader will find ``linearly
convex,'' ``lineally convex,'' ``planarkonvex,'' and ``complex geometric convexity'' used to describe the
same or similar notions.)  

In [Wh] we proved Helton and Marshall's conjecture with some
unnatural assumptions: that $K$ can be smoothly deformed through sets
$K^t\subset\Gamma\times\C^n$ such that for all $t$, $K^t_z$ is strictly hypoconvex; that for
$t$ large, $K^t_z$ is a ball of radius $t$ and for $t$ small, $K^t_z$
shrinks around the origin in $\C^n$.  Further, we assumed that the
polynomial hull of $K$ contains at least one analytic graph.  Here we
show that these assumptions are unnecessary by
use of some constructions in work of Lempert [L1,L2,L3,L4,L5,L6]; we will show that the required deformation can be constructed.
This is the subject of Theorem \hullthm.  Let
$\Pi:\Delta\times\C^n\rightarrow\Delta$ be projection.  We assume the following about $K$:

$$\eqalign{(1) & \left[\vbox{\hsize=5 in $K$ is a $C^\infty$ submanifold of
$\Gamma\times\C^n$ parametrized by a $C^\infty$ diffeomorphism $I:\Gamma\times\partial B_n\rightarrow K$, $I(z,w)=(z,I_2(z,w))$ for some $I_2$, such that for every $z\in \Gamma$, $K_z$ is the boundary of a bounded connected strictly hypoconvex
open set $Y_z$ in $\C^n$.
\vskip -0.2 in
}\right.}$$

{\bf Theorem \hullthm.} {\it Suppose $K$ satisfies (1). Then $\widehat K\cap
\Pi^{-1}(\opdisk)$ is either (a) empty, (b) the graph of a single analytic
mapping $f:\opdisk\rightarrow\C^n$ which extends to be in $C^\infty(\Delta)$ or (c) the union of infinitely many
analytic graphs over the disk which extend to be in $C^\infty(\Delta)$.  If $(z_0,w_0)\in\partial\widehat K,$
$|z_0|<1$, there is only one analytic $f\in H^\infty(\Delta)^n$ such that $f(z_0)=w_0$ and $f(z)\in \widehat K_z$ for a.{}e.{} $z\in\Gamma$.  In case (c), we
have that for all $z\in\Delta $, $\widehat K_z$ is hypoconvex with
$C^1$ boundary.}

The smoothness statement in the last sentence can be improved; 
see [BL,Wh].

We use similar analysis to study another problem posed in [HMa].  Let
$\rho:\Gamma\times\C^n\rightarrow \R$.  A problem of $H^\infty$
(frequency domain) design is to find $f\in H^\infty(\Delta)^n$ such
that the essential supremum of $z\mapsto\rho(z,f(z))$ is as small as
possible.  (In $H^\infty$ design, $f$ determines a design for a
machine and $\rho$ measures something bad about that design that we
wish to keep small.  See [HMa] for a discussion of this.)  The
question we address here is whether a minimizer $f$ is unique and
whether it is smooth up to the boundary of $\Delta$.  Let
$$\gamma_\rho\equiv\inf_{f\in H^\infty(\Delta)^n}\hbox{\rm
ess}\sup_{z\in\Gamma}\rho(z,f(z)).\tag 2$$

{\bf Theorem \optthm .} {\it Suppose that
$\rho:\Gamma\times\C^n\rightarrow\R$ is continuous and that the set
$K$ where $\rho $ equals $\gamma_\rho$ satisfies (1), where 
$Y_z=\{w\in\C^n:\rho(z,w)<\gamma_\rho\}$.  Then there exists a unique
$\phi\in H^\infty(\Delta )^n$ such that $\rho(z,\phi(z))\leq
\gamma_\rho$ for a.e. $z\in\Gamma$; in fact $\phi$ extends to be in $C^\infty(\Delta )$ and $\rho(z,\phi(z))=\gamma_\rho$ for all
$z\in\Gamma$.  The set $\widehat K\cap\Pi^{-1}(\opdisk)$ contains only
the points on the graph of $\phi$.  If $\rho(\overline z,\overline
w)=\rho(z,w)$, then $\phi $ is
\R$^n$-valued on the real axis, i.e. $\overline {\phi(\overline
z)}=\phi(z)$.}

This generalizes our own work [Wh, Theorem 3] and work of Vityaev [V,
Theorem 1.5], which assumes that an optimizer $\phi$ is smooth up to
the boundary of $\Delta$ and proves that $\phi$ is unique in the class
of analytic mappings $\opdisk\rightarrow\C^n$ which extend to be
smooth on the boundary.

We shall first prove a slightly weaker version of Theorem \hullthm:

{\bf Theorem \fathull.} {\it Suppose $K$ satisfies (1) such that for every $z\in \Gamma$, $Y_z$  contains the origin.  Then $\widehat
K\cap
\Pi^{-1}(\opdisk)$ is the union of graphs of analytic mappings
$f:\opdisk\rightarrow
\C^n$ which extend to be in $C^\infty(\Delta)$.  If $(z_0,w_0)\in\partial\widehat K,$
$|z_0|<1$, there is only one analytic $f\in H^\infty(\Delta)^n$ such that $f(z_0)=w_0$ and $f(z)\in \widehat K_z$ for a.{}e.{} $z\in\Gamma$.  For all $z\in\Delta $, $\widehat
K_z$ is hypoconvex with $C^1$ boundary.}

{\bf Remark.} Theorem \fathull is true if we replace the requirement that $K$
enclose the zero graph by the requirement that it enclose some
analytic graph, by holomorphic change of coordinates.  Theorem
\hullthm and Theorem \optthm are known in the case $n=1$, so we assume
throughout that $n>1$.  We first prove Theorem \fathull, then Theorem
\optthm , and then Theorem \hullthm.

We discuss some notation.
Let $Df(x)$ be the gradient of a function $f$ at $x$ and
let $D^2f(x)$ be the real Hessian (i.e., the second derivative) at
$x$.  We write $D^2f(x)[u,v]$ to denote the value of the bilinear
form $D^2f(x)$ at the pair $u,v$.  Let $\rho $ be a real valued
function on an open $U\subset\Gamma\times\C^n$.  Then $D_w\rho(z,w)$
is the vector $({\partial\rho
\over\partial w_1}(z,w),{\partial\rho
\over\partial w_2}(z,w),...,{\partial\rho
\over\partial w_n}(z,w))$.  We say that $\rho
:U\rightarrow
\R$ is strictly hypoconvex on $U$ if $D_w\rho\neq 0$ there and for
every compact subset of $U$ there exists $\kappa>0$ such that
whenever
$$\sum_{j=1}^n u_j {\partial\rho\over \partial w_j}
(z,w_1,w_2,...,w_n)=0$$ on that compact subset we have
$$D^2\rho(z,w_1,w_2,...,w_n)[0,u_1,u_2,...,u_n][0,u_1,u_2,...,u_n]\geq
\kappa|u|^2.$$

In Theorem \hullthm and Theorem \optthm , the fibers $K_z$ of $K$
enclose bounded connected strictly hypoconvex sets which are also
polynomially convex; see [Wh].
Work of Lempert shows that the compact set bounded by $K_z$ is
homeomorphic to $\overline B_n$ via a mapping $:\overline B_n\rightarrow\widehat K_z$ which is $C^\infty$ except at the origin.
(See also
Corollary 4.6.9 and Theorem 4.6.12 of [H\" o].)  Since the closed set
bounded by $K_z$ is polynomially convex for $z\in\Gamma$, when we
compute the polynomial hull of $K$, its fiber over $z\in\Gamma$ will
be $(\widehat K)_z=\widehat{K_z}$ by Lemma 1 of [HMa].  For ease in
notation, we shall refer to this set as $\widehat K_z$.

We say that $K$ is a level set of $\rho$ if there exists a constant $c$ such that $K$ is the set of all points in the domain of $\rho$ where $\rho $ takes the value $c$.  We shall
frequently need a defining function for a set satisfying (1).  Assuming that ${\Cal K}$ satisfies (1) and ${\Cal Y}_z$ is the bounded open set bounded by ${\Cal K}_z$, here are properties
of $\rho$ that we shall need:

$$\eqalign {(3)&\left [\vbox {\hsize =4.5 in\vbox{(a) $\rho $ is
defined and $C^\infty$ in a neighborhood of ${\Cal K}$ in $\Gamma\times\C^n$, such that $D_w\rho(z,w)\neq 0$, ${\Cal K}$ is a level set of $\rho$, 
and $\rho$ is strictly hypoconvex on its domain;} 
\vbox{(b) The values of $\rho(z,\cdot)$
are (strictly) smaller in ${\Cal Y}_z$ than on ${\Cal K}_z$.}
\vskip -0.35in
}
\right .}$$ 

For our original $K$, we may find such a $\rho $ in the following manner: let $d_z(w)$ be
the distance from $w$ to $K_z$.  Then for $(z,w)$ near $K$, let $\rho(z,w)$ be $1+d_z(w)$
for $w$ in the unbounded component of $(K_z)^c$ and $1-d_z(w)$ for $w$
in the bounded component of $(K_z)^c$.  This $\rho$ is $C^\infty$ and strictly hypoconvex for $(z,w)$ sufficiently close to $K$; we can then extend $\rho$ to
$\Gamma\times\C^n$, so that $\rho $ satisfies (3).

We seek a $\rho$ which satisfies

$$\eqalign {(4)&\left [\vbox {\hsize =4.5 in\vbox{(a)
$\rho:\Gamma\times\C^n\rightarrow[0,\infty)$ is continuous, and
$C^\infty$-smooth where $\rho\neq 0$;}\medskip
\vbox {(b) there exists an $R>1$ such that if $0<\rho(z,w)\leq R$
then $D_w\rho(z,w)\neq 0$ and $\rho$ is strictly hypoconvex as 
in (3) where $\rho $ is smooth;}\medskip
\vbox{(c) for every $t$, $0<t\leq R$, the set $K^t$ where $\rho=t$ is
compact, with fibers $K^t_z$ diffeomorphic to $\partial
B_n$, where $B_n$ is the open unit ball in $\C^n$;}
\vbox{(d) $K=K^1$;}
\vbox{(e) $\{(z,w) \bigl| \rho(z,w)=R\}=\{(z,w) \bigl | |w|=R\}$ and 
$\rho(z,w)>R$ if $|w|>R$;}\medskip
\vbox{(f) There exists a continuous function $S(z)$, $|S(z)|<R$, 
such that $\{(z,w)\in
\Gamma\times\C^n\bigl| w=S(z)\}=\{(z,w)\in
\Gamma\times\C^n\bigl| \rho(z,w)=0\}$.}
\vskip -0.9 in
}\right .}$$

(Such a $\rho$ will satisfy (3) on p.{}679 of [Wh].)

$\S 2$ {\bf Lempert's work.}

Let $A(\Delta)=\{f:\Delta\rightarrow\C |
\hbox{$f$ is continuous on $\Delta$ and analytic in \opdisk}\}$.
We shall rely on constructions contained in work of Lempert, which we
shall summarize.  (See [L1,L2,L3,L4,L5,L6].)  This work will
allow us to deform the fibers $K_z$ of our set $K$ in
a manner which is necessary to obtain a defining function satisfying (4).

Fix any $z\in\Gamma $.  As observed earlier, the interior of $\widehat K_z$ is a
closed $C^\infty$-bounded strictly hypoconvex domain $Y_z$.  Let $\monge
:\widehat K_z\rightarrow\R$ be Lempert's solution to the homogeneous
complex Monge-Amp\` ere equation for $\widehat K_z$ with logarithmic
singularity at $S(z)$ ($\monge (w)=\log|w-S(z)|+O(1)$ as
$w\rightarrow S(z)$) and boundary value $0$ on $K_z$.  (See
[L2,L3].)  Let $\ui(z,w)\equiv (\ui)_z(w)=e^{2\monge (w)}$.  According to Lempert's
work, given any pair $(p,\nu)$ with $p$ in int $\widehat K_z$ and $\nu\in\C^n$, $|\nu|=1$, there
exists a unique extremal analytic disc, some of whose properties we
shall need (see [L2]).  An extremal disc satisfies the following (among other things): it is a continuous $f_{z,\nu}:\Delta\rightarrow
\widehat K_z$ such that $f_{z,\nu}$ is analytic on \opdisk, $f_{z,\nu}(0)=p$, $f'_{z,\nu}(0)$ is a positive multiple of $\nu$, $f_{z,\nu}(\partial\Delta)\subset
K_z$ and $f_{z,\nu}$ has a holomorphic left inverse $F_{z,\nu}$ which is $C^\infty$ on $\widehat K_z$ (see Proposition 1 of [L2] and its
proof).  We will let $p=S(z)$.  The level sets of $F_{z,\nu}$ 
are complex hyperplanes (intersected
with $\widehat K_z$).  Each of these complex hyperplanes meets the
extremal disc in exactly one point([L2], Proposition 1 and its proof).
At this point, the complex hyperplane is tangent to the level set of
$(\ui)_z$ through that point.  On the extremal disc parametrized by $f_{z,\nu}$,
$\log|F_{z,\nu}(w)|=\tau _z(w)$, so $|F_{z,\nu}(w)|^2=e^{2\tau _z(w)}=\ui (z,w)$.
Furthermore, on $\widehat K_z$, $\tau _z(w)=\max_{\nu\in C^n, |\nu|=1}\log|F_{z,\nu}(w)|$.  For $0<c\leq 1$ the set
$\{w\in\widehat K_z | \tau _z(w)\leq c\}$ is homeomorphic to the closed
ball via a mapping which is $C^\infty$ except at $p$; Lempert constructs this homeomorphism in [L3] for convex
domains.
It is also true (see the Reflection Principle in [L2, p.{}348]) that the 
extremal disks for $K_z$
such that $f(0)=S(z)$ extend to be in $C^\infty(\Delta)$.  For $0<\ui\leq 1$, $\ui$ is strictly hypoconvex.
([L4,pp.{}523,576])

We also make use of Lempert's solution $\tilde\tau_z$ to the homogeneous complex
Monge-Amp\` ere equation on the unbounded component of $(K_z)^c$,
(with boundary value $0$ on $K_z$).  (See [L5, p.{}881] and the remark
extending the result to hypoconvex domains on p.{}884.)  Let $\uo(z,w)=
e^{2\tilde\tau_z(w)}$.  Then $\uo$
extends continuously to $K$, attaining the boundary value $1$.
([L5,pp.{}881,884])  Where $\uo>
1$, $\uo $ is strictly hypoconvex.([L5,p.{}883])

$\S 3$ {\bf Deformation of the set $K$.}

{\bf Lemma \lemcenter.} {\it If $K$ satisfies (1), there exists a $C^\infty$ mapping $S:\Gamma\rightarrow\C^n$ such that $S(z)\in Y_z$ for all $z\in\Gamma$.}

{\it Proof.}  This follows almost immediately from the definition of $K$.  Fix $q\in\partial B_n$.  There we have that $I_2(z,q)\in K_z$ for all $z\in\Gamma$; let $n(z,w)$ denote the inward pointing unit normal to $K_z$ at $w\in K_z$.  Since $I_2(z,q)$ is $C^\infty$ in $z$, we may let $S(z)=I_2(z,q)+r n(z,I_2(z,q))$ if $r$ is chosen sufficiently small.  (This `pushes' $I_2(z,q)$ a short distance into $Y_z$.) $\square$

We assume $S(z)$ is any function satisfying the conditions of Lemma \lemcenter.
For $\epsilon>0$, fixed $z\in\Gamma$, $\nu\in \C^n$, let $f_{z,\nu}:\Delta\rightarrow \widehat
K_z$ be the  extremal mapping with $f(0)=S(z)$ and $f'(0)=\lambda \nu$, $\lambda>0$.  Then let $F_{z,\nu}$ be the holomorphic left inverse to $f_{z,\nu}$ as defined in [L2],
and let $F_{z,\nu}^\epsilon(w)^2= |F_{z,\nu}(w)|^2+\epsilon|w-f_{z,\nu}(F_{z,\nu}(w))|^2$.

{\bf Lemma \lemsup.}  {\it There exists $\epsilon >0$ such that for all $z\in\Gamma$, on the set $\{w\in\widehat K_z: (\ui)_z(w)\leq 1/2\}$ we have that
$(\ui)_z(w)$ is the maximum
$$\max_{\nu\in C^n, |\nu|=1}
(F_{z,\nu}^\epsilon)(w)^2.$$  (The set $\{w\in\widehat K_z: (\ui)_z(w)\leq 1/2\}\subset\subset$ int $\widehat
K_z$.)}

{\it Proof.}  The maximum above exists because the quantity there depends continuously on $\nu$, which varies over a compact set.  Let $\nu\in\C^n$, $|\nu|=1$.  If $w$ is in the image of $f_{z,\nu}$, then
$F_{z,\nu}^\epsilon(w)^2 =|F_{z,\nu}(w)|^2$ since $w=f_{z,\nu}(F_{z,\nu}(w))$.  The function $(\ui)_z(w)$ is less than or equal to $\max_{\nu\in C^n, |\nu|=1}
(F_{z,\nu}^\epsilon)(w)^2$ because Lempert proves that
$\tau _z(w)=\max_{\nu\in C^n, |\nu|=1}\log|F_{z,\nu}(w)|$, so $(\ui)_z(w)=\max_{\nu\in C^n,|\nu|=1}|F_{z,\nu}(w)|^2\leq \max_{\nu\in C^n,|\nu|=1}
(F_{z,\nu}^\epsilon)(w)^2$.  To prove the reverse inequality, we
choose $\epsilon $ so small that $2\epsilon $ is less than $D^2\ui
(z,w)[v,v]$, where $z\in\Gamma$,
$1/3\leq (\ui)_z(w)\leq 2/3$, $|v|=1$ and $v$ is a complex tangent at $w$ to the
surface $\{w:(\ui)_z(w)= c\}$, for some $c$, $1/3\leq c\leq 2/3$.  (This is possible since $\ui $
is $C^\infty$ and strictly hypoconvex where $0<\ui<1$.)  
Then for each $f_{z,\nu}$ with
holomorphic left inverse $F_{z,\nu}$, we claim $|F_{z,\nu}(w)|^2+\epsilon|w-f_{z,\nu}(F_{z,\nu}(w))|^2\leq (\ui)_z(w)$ in a neighborhood of $f(\{\lambda\in\C||\lambda|=1/\sqrt 2\})\subset \{w\in\widehat K_z:(\ui)_z(w)=1/2\}$.  To see why the claim holds, first note that if $w$ is in the image of $f_{z,\nu}$, then $w=f_{z,\nu}(\lambda)$ for some $\lambda$ in the open unit disk; substituting $w=f_{z,\nu}(\lambda)$ into $\epsilon|w-f_{z,\nu}(F_{z,\nu}(w))|^2$ yields zero.  Then from Lempert's work we already know that $|F_{z,\nu}(w)|^2= (\ui)_z(w)$ for $w$ in the image of $f$.  Thus the claim holds - in fact, is an equality - for $w$ in the image of $f_{z,\nu}$.  Now let $w$ be a point of the form $f_{z,\nu}(\lambda)$ for some $\lambda$ such that $\sqrt{1\over 3}\leq |\lambda|\leq\sqrt{2\over 3}$.  Then the second
derivative of $\ui (z,\cdot)=\exp(2\tau _z(\cdot))$ at $w$ in a
direction $v$ along a level set of $F_{z,\nu}$ is $> 2\epsilon $.  The first derivatives of $\ui(z,\cdot)$ and $F_{z,\nu}(w)^2$ at $w$ in the direction of $v$ are both zero.  (We recall that a level set of $F_{z,\nu}$ is
an affine complex hyperplane tangent to a level set of $(\ui)_z$ at a point in the image of $f_{z,\nu}$.
Every point in the image of $f_{z,\nu}$ lies on precisely one of these
level sets of $F_{z,\nu}$.)  By comparing derivatives of functions, elementary calculus allows us to conclude that $|F_{z,\nu}(x)|^2+\epsilon|x-f_{z,\nu}(F_{z,\nu}(x))|^2< (\ui)_z(x)$ if $\sqrt{1\over 3}\leq |F_{z,\nu}(x)|\leq\sqrt{2\over 3}$ and $x$ is sufficiently close to the image of $f_{z,\nu}$.  We may then choose a $\delta>0$ independent of $\nu$ and $z$
such that for all $(z,\nu)$,
$F_{z,\nu}^\epsilon(w)^2\leq (\ui)_z(w)$ on a $\delta$-neighborhood of
$f_{z,\nu}(\{\lambda\in\C||\lambda| =1/\sqrt 2\})$; the above argument achieves this in a neighborhood of a fixed $(z,\nu)$.  Then the continuity of the functions involved and a compactness argument obtains the $\delta$ for all $(z,\nu)\in\Gamma\times\partial B_n$.  We may further choose $\epsilon $ even
smaller that $$F_{z,\nu}^\epsilon(w)^2\leq (\ui)_z(w)\tag 5$$ for any $w$
{\it outside} the $\delta$ - neighborhood determined above independent of $\nu$
and $z$ such that $\ui (z,w)=1/2$.  (We may do this for any fixed $(z,\nu)$ because on the topological sphere where $(\ui)_z(w)=1/2$, $|F_{z,\nu}(w)|^2<(\ui)_z(w)$ outside a $\delta$-neighborhood of $f_{z,\nu}(\{\lambda\in\C||\lambda| =1/\sqrt 2\})$, so the $\epsilon$ can be found for this fixed $(z,\nu)$.  This $\epsilon$ will work for $(z',\nu')$ nearby, by continuity.  Since $(z,\nu)$ is allowed to vary in the compact set $\Gamma\times\partial B_n$, the $\epsilon $ can be found so that (5) holds, as desired.)  Combining the two neighborhoods, we find that
on the set where $(\ui)_z(w)= 1/2$,
$F_{z,\nu}^\epsilon(w)^2\leq (\ui)_z(w)$.  We claim that this
inequality extends to the region where $\ui \leq 1/2$.  We have
$|F_{z,\nu}(w)|^2+\epsilon |w-f_{z,\nu}(F_{z,\nu}(w))|^2\leq (\ui)_z(w)$ if $(\ui)_z(w)=1/2$.  By Lempert's work, if $w\in Y_z\equiv int\,\widehat K_z$, there is a $\mu\in\C^n$, $|\mu|=1$ such that $w=f_{z,\mu}(\lambda)$, $|\lambda|<1$, for some extremal
$f_{z,\mu}$; $F_{z,\mu}$ is the holomorphic left inverse of $f_{z,\mu}$.  Then
$(\ui)_z(f_{z,\mu}(\lambda))=|\lambda|^2$, since if $w=
f_{z,\mu}(\lambda)$, $\tau _z(f_{z,\mu}(\lambda))=\log|F_{z,\mu}(
f_{z,\mu}(\lambda))|=\log|\lambda|$.  Also, since $F_{z,\nu}(f_{z,\mu}(0))=0$ and
$f_{z,\mu}(0) -f_{z,\nu}(F_{z,\nu}(f_{z,\mu}(0)))=\origin$, we can write that for
$|\lambda|=1/\sqrt 2$,
$|\lambda|^2|A(\lambda)|^2+|\lambda|^2|B(\lambda)|^2\leq|\lambda|^2$,
where $A,B$ are analytic on the disk, $F_{z,\nu}(f_{z,\mu}(\lambda)) =\lambda
A(\lambda)$, and $f_{z,\mu}(\lambda) -f_{z,\nu}(F_{z,\nu}(f_{z,\mu}(\lambda))) =\lambda
B(\lambda)$.  Then $|A(\lambda)|^2+|B(\lambda)|^2\leq 1$ for
$|\lambda|=1/\sqrt 2$.  Since the left side is subharmonic in
$\lambda$, this inequality extends to the disk where $|\lambda|\leq
1/\sqrt 2$.  Multiplying back by the $|\lambda|^2$, we obtain the
desired inequality on the image of $f_{z,\mu}$.  Since every point $w$
such that $\ui (z,w)\leq 1/2$ lies on such an extremal disk, we find
that the claim holds: $F_{z,\nu}^\epsilon(w)^2\leq (\ui)_z(w)$ if
$\ui (z,w)\leq 1/2$.  We let $(\tilde{\ui})_z(w)=
\max_{\nu\in C^n,|\nu|=1} F_{z,\nu}^\epsilon(w)^2$.  Thus what we have just done shows that
$(\tilde{\ui})_z(w)\leq (\ui)_z(w)$ if $\ui (z,w)\leq 1/2$.
As noted previously, $(\ui)_z(w)\leq (\tilde{\ui})_z(w)$ on
$\widehat K_z$.  Combining these, we have the desired result.
$\square$

{\bf Lemma \lemdelta.}  {\it There exists a $\delta>0$ such that $\ui(z,w)$ is strictly convex in $w$ for $0<|w-S(z)|<\delta$ and all $z$; hence
$\{(z,w)\in\Gamma\times\C^n\bigl| \ui (z,w)\leq
\delta\}$ has strictly convex fibers over the circle.}

{\it Proof.}  Let $f_{z,\nu}$ be the extremal disc for $\widehat K_z$ such that
$f'(0)$ is a positive multiple of $\nu \in\C^n$ (see [L2, Theorem 3])
and $f(0)=S(z)$.  Let $F_{z,\nu}$ be the holomorphic left inverse of
$f_{z,\nu}$ and let
$F^\epsilon_{z,\nu}(w)^2=|F_{z,\nu}(w)|^2+\epsilon|w-f_{z,\nu}(F_{z,\nu}
(w))|^2$, where $\epsilon$ is as chosen in Lemma \lemsup.  We claim
that there exists $\delta>0$ independent of $z,\nu $ such that
$(F^\epsilon_{z,\nu})^2$ is strictly convex on $B_n(S(z),\delta)$,
the open ball of radius $\delta$ centered at $S(z)$.  If not,
choose $\delta_n\rightarrow 0$ and a corresponding sequence of
$(F^{\epsilon}_{z_n,\nu_n})^2$ such that the claim is false.  Passing to
a limit, we find that there exists a pair $(z,\nu)$ such that
$(F^\epsilon_{z,\nu}(w))^2$ is not strictly convex in $w$ at $S(z)$.  We show that this is impossible.  We have that (a) the real
Hessian of $|F_{z,\nu}|^2$ is nonnegative at the origin (since
$|F_{z,\nu}|^2$ has a local minimum there) and positive in every
direction except those complex directions orthogonal to $\nu$ (recall
that $D_wF_{z,\nu}$ is nonzero, as the left inverse of $f_{z,\nu}$.)
In the complex directions orthogonal to $\nu$, $F$ is identically
zero.  Further, (b) $\epsilon|w-f_{z,\nu}(F_{z,\nu}(w))|^2$ is
strictly convex in those directions complex orthogonal to $\nu$.  Thus
the sum $|F_{z,\nu}(w)|^2+\epsilon|w-f_{z,\nu}(F_{z,\nu}(w))|^2$ is
strictly convex at $S(z)$, a contradiction.  Thus $\ui (z,\cdot)$ is strictly convex in
$w$ for sufficiently small $|w-S(z)|\neq 0$ as well, as a maximum of the
functions $(F^\epsilon_{z,\nu})^2$ (Lemma \lemsup).  $\square$


In an earlier version of this paper, we noted that if a strictly hypoconvex domain were completely circled about the origin, then its Kobayashi balls about the origin would then be dilations of the domain.  Then Lemma \lemdelta shows that a strictly hypoconvex domain which is completely circled must be convex.  We indicated then that there must be a more elementary proof of this fact.  Miran \v Cerne 
has provided precisely such a proof in [\v Ce2].

We will now piece together various defining functions for $K$ in order to obtain one which satisfies (4); we will need the following lemma regarding mollification of such functions.

{\bf Lemma \lemmollify .}  {\it Suppose that $K$ satisfies (1), has defining functions $\rho_1,\rho_2$ which satisfy (3), $\rho_1=\rho_2=1$ on $K$, and $|D_w\rho_1|>|D_w\rho_2|$ on $K$.  Let $\rho(z,w)=\rho_1(z,w)$ if $\rho_1(z,w)>1$, and $\rho(z,w)=\rho_2(z,w)$ if $\rho_1(z,w)\leq 1$.  Let $\phi_\epsilon$ be a $C^\infty$ approximation to the convolution identity in $(z,w)\in\Gamma\times\C^n$ whose support has radius $<\epsilon$.  Let $N(K,\epsilon)=\{(z,w)\in\Gamma\times\C^n:|\rho(z,w)-1|\leq\epsilon\}$ and let $K_\epsilon^c\equiv\{(z,w)\in\Gamma\times\C^n:(\rho*\phi_\epsilon)(z,w)=c\}$.  Then given any neighborhood $M$ of $K$, there exists a neighborhood $N(K,\epsilon_1)$ of $K$ contained in $M$ and $\epsilon_2 $ sufficiently small such that $N(K,\epsilon_1/4)\subset \{(z,w)\in\Gamma\times\C^n:|\rho*\phi_\epsilon(z,w)-1|<\epsilon_1/2\}\subset N(K,\epsilon_1)$ for $\epsilon\leq\epsilon_2$, and $K^\epsilon_c$ satisfies (1) for $c\in[1-\epsilon_1/2,1+\epsilon_1/2]$.  Furthermore, $\rho*\phi_\epsilon$ is strictly hypoconvex on the set $\{(z,w)\in\Gamma\times\C^n:|\rho*\phi_\epsilon(z,w)-1|<\epsilon_1/2\}$.}

{\bf Note.}  An analogous result holds in the event that we only have $\rho_1=\rho_2=k$ on $K$ for some real constant $k$; merely apply Lemma \lemmollify to $\rho_1-k+1,\rho_2-k+1$, which both have value $1$ on $K$.  

{\it Proof.}  Note that the gradient condition on $\rho_1,\rho_2$ guarantees that $\rho$ is the maximum of $\rho_1$ and $\rho_2$ in some neighborhood $N(K,2\epsilon_1)\subset M$ of $K$.  We will convolve $\rho$ with a close approximation $\phi_\epsilon(z,w)$ to the identity in $(z,w)\in\Gamma\times\C^n$, where $\epsilon>0$.  We assume $\epsilon_2$ is small enough that $\rho$ does not vary more than $\epsilon_1/4$ on a disk of radius $\epsilon_2$ which meets $N\equiv N(K,\epsilon_1)$.  Since $\rho*\phi_\epsilon$ converges locally uniformly to $\rho$ as $\epsilon\rightarrow 0$, we obtain the statement that for $\epsilon_2$ small, $N(K,\epsilon_1/4)\subset \{(z,w)\in\Gamma\times\C^n:|\rho*\phi_\epsilon(z,w)-1|<\epsilon_1/2\}\subset N(K,\epsilon_1)$ for $\epsilon\leq\epsilon_2$.  Now $D_w\rho_1$ and $D_w\rho_2$ are nonzero in $N$ and have the same direction on $K$, so the unit normal to the level sets of $\rho$ varies continuously near $N$.  Assume $\epsilon_2$ is small enough that on a ball of radius $\epsilon_2$ which meets $N$, the angle of the unit normal to level sets of $\rho(z,\cdot)$ varies no more than $\pi/10$.   Then $D_w(\rho*\phi_\epsilon)$ is nonzero near $N$ for $\epsilon\leq\epsilon_2$.  Let $T(z,w)$ be the complex tangent space in $\C^n$ to the level set of $\rho(z,\cdot)$ at $w$.  Then $T$ is continuous.  By choosing $\epsilon$ even smaller, we may make the direction of ${D_w(\rho*\phi_\epsilon)\over |D_w(\rho*\phi_\epsilon)|}$ as (uniformly) near to the directions ${D_w\rho_1\over |D_w\rho_1|}$ (on $\{\rho_1\geq 1\}\cap N$) and ${D_w\rho_2\over |D_w\rho_2|}$ (on $\{\rho_1\leq 1\}\cap N$) as desired.  This implies that for $\epsilon$ small, level sets of $\rho*\phi_\epsilon$ which meet $N$ have spherical fibers over $\Gamma$.  For $\epsilon$ sufficiently small and $1-\epsilon_1/2\leq c\leq 1+\epsilon_1/2$, $K^c_\epsilon$ may be parametrized by a diffeomorphism from $\Gamma\times\partial B_n$ by composing the map which parametrizes $\{(z,w):\rho(z,w)=c\}$ by orthogonal projection along its normals to $K^c_\epsilon$.  Now $\rho_1(z,\cdot),\rho_2(z,\cdot)$ are strictly convex on $T(z,w)$ at $w$.  In fact, for all $(z,w)\in N(K,\epsilon_1)$ there exists $\theta>0$ such that they are strictly convex at $w$ in any direction with an angle within $\theta$ of $T(z,w)$.  Assume $\epsilon_2$ so small that on a ball of radius $\epsilon_2$ meeting $N$, the normal to $T(z,w)$ varies no more than $\theta/3$.  Assume $\epsilon_2$ so small that for $\epsilon\leq\epsilon_2$, the direction of $D_w(\rho*\phi_\epsilon)$ is within $\theta/3$ of $ {D_w\rho_1\over |D_w\rho_1|} $ on $\{\rho_1\geq 1\}\cap N$ and ${D_w\rho_2\over |D_w\rho_2|} $ on $\{\rho_1\leq 1\}\cap N$.  Since $\rho_1,\rho_2$ are $C^\infty$ and strictly hypoconvex on $N$, this means that for $\epsilon_2$ small enough, there exists $C>0$ such that $${\rho_i(z,w+hv)+\rho_i(z,w-hv)\over 2}\geq \rho_i(z,w)+Ch^2,\tag 6$$ for all $h\in\R$, $|h|\leq\epsilon_2$, $(z,w)$ near  $N$, $i=1,2$ and $v$ a unit vector in $\C^n$ whose direction is within $\theta$ of $T(z,w)$.  The property (6) is preserved under taking the maximum of $\rho_1,\rho_2$, so $\rho$ satisfies $${\rho(z,w+hv)+\rho(z,w-hv)\over 2}\geq \rho(z,w)+Ch^2,\tag 7$$ for all $h\in\R$, $|h|\leq\epsilon_2$, $(z,w)$ near $N$, $i=1,2$ and $v$ a unit vector in $\C^n$ whose direction is within $\theta$ of $T(z,w)$.  For $(z,w)$ near $N$, if $v$ is within $\theta/3$ of being a tangent to the level set of $\rho*\phi_\epsilon(z,\cdot)$ at $w$, then it is within $2\theta/3$ of $T(z,w)$, so within $\theta$ of $T(z,w+h\tilde v)$ if $|\tilde v|=1$, $|h|\leq \epsilon_2$.  Thus for $(z_0,w_0)$ fixed near $N$, $${\rho(z,w+hv)+\rho(z,w-hv)\over 2}\geq \rho(z,w)+Ch^2,\tag 8$$ for all $h\in\R$, $|h|\leq\epsilon_2$, $(z,w)$ within $\epsilon_2$ of $(z_0,w_0)$, $i=1,2$ and $v$ a unit vector in $\C^n$ whose direction is within $\theta/3$ of the tangent space to the level set of $\rho*\phi_\epsilon(z_0,\cdot)$ at $w_0$.  Then when convolving $\rho$ with $\phi_\epsilon$, $\epsilon\leq\epsilon_2$, $\rho*\phi_\epsilon$ satisfies $${\rho*\phi_\epsilon(z,w+hv)+\rho*\phi_\epsilon(z,w-hv)\over 2}\geq \rho*\phi_\epsilon(z,w)+Ch^2,\tag 9$$ for all $h\in\R$, $|h|\leq\epsilon_2$, $(z,w)$ near $N$, $i=1,2$ and $v$ a unit vector in $\C^n$ at $w$ whose direction is within $\theta/3$ of being tangent at $w$ to the level set of $\rho*\phi_\epsilon(z,\cdot)$.  This shows that $\rho*\phi_\epsilon$ is strictly hypoconvex in $N$.  Since the level sets of $\rho*\phi_\epsilon$ were already shown to have spherical fibers over $\Gamma$, they now must enclose bounded strictly hypoconvex domains, so must satisfy (1), as desired.  $\square$

{\bf Lemma \lemlocmollify .} {\it Let $K,\rho_1,\rho_2,\rho$ be as defined in Lemma \lemmollify.  Then given a neighborhood $M$ of $K$ there exists a strictly hypoconvex function $\rho_3$ with the same domain as $\rho$ whose level sets satisfy (1) and which coincides with $\rho$ outside of $M$.}

{\bf Note.}  Lemma \lemlocmollify stills holds in the event that we only have $\rho_1=\rho_2=k$ on $K$ for some real constant $k$; simply apply the lemma to the functions $\rho_1-k+1,\rho_2-k+1$.

{\it Proof.}  Apply Lemma \lemmollify to obtain $N=N(K,\epsilon_1)$ and $\epsilon_2$.  Using a partition of unity, we may write $\rho(z,w)$ as a sum of two functions $\Phi,\Psi$, one of which (say $\Phi$) has its support in $N(K,4\epsilon_1/5)=\{(z,w)\in\Gamma\times\C^n:|\rho(z,w)-1|\leq4\epsilon_1/5\}$, and $\Phi$ is identically $\rho$ in a neighborhood of $N(K,3\epsilon_1/5)$.  If $\epsilon_2$ is small enough, we find $\Phi*\phi_\epsilon=\rho*\phi_\epsilon$ in $N(K,\epsilon_1/2)$, so is strictly hypoconvex there if $\epsilon\leq\epsilon_2$, as indicated in Lemma \lemmollify.  Furthermore, $K^c_\epsilon$ satisfies (1) for $c\in[1-\epsilon_1/2,1+\epsilon_1/2]$.  But if $\epsilon_2$ is even smaller, $\Phi_\epsilon+\Psi$ is strictly hypoconvex outside of $N(K,\epsilon_1/4)$ for $\epsilon\leq\epsilon_2$ because $\rho$ is strictly hypoconvex outside of $N(K,\epsilon_1/4)$, its level sets satisfy (1), and $\Phi_\epsilon\rightarrow\Phi$ in local $C^2$ norm.  Then $\rho_3\equiv\Phi_\epsilon+\Psi$ satisfies the conditions of the lemma for $\epsilon\leq\epsilon_2$. $\square$

Let us recall the function $S(z)$ whose existence was asserted in Lemma \lemcenter.

{\bf Lemma \lemdeffuncA .}  {\it Let $K$ satisfy (1).  Then there exists a strictly hypoconvex defining function $\uo$ defined in a neighborhood of $\{(z,w)\in\Gamma\times\C^n:w\in\widehat K_z\}$ which is strictly hypoconvex where it is nonzero, and such that for $t$ sufficiently small, $\{w\in\C^n:\uo(z,w)=t\}$ is a sphere of radius independent of $z\in\Gamma$ centered at $S(z)$.  The zero set of $\uo(z,\cdot)$ on $\widehat K_z$ consists only of $S(z)$, and $\{(z,w)\in\C^n:\uo(z,w)=t\}$ satisfies (1) for $t>0$.}

{\it Proof.}  Lempert [L2, Prop. 13] proves the existence of a $C^\infty$ mapping $\Xi:\Gamma\times\C^n\rightarrow A(\Delta)^n$ such that $\Xi(z,\nu)$ is the extremal disc $\phi$ for $\widehat K_z$, where $\phi(0)=S(z)$ and $\phi'(0) $ is a positive multiple of $\nu$.  Composing with evaluation at the real number $r\in\Delta$ (denoted $\Xi(z,\nu)(r)$), $(z,\nu)\mapsto (z,\Xi(z,\nu)(r))$ becomes a $C^\infty$ parametrization of the set $\{(z,w)\in\Gamma\times\C^n:u_1(z,w)=r^2\}$.  (See \S 2 for the definition of $\ui$.)  Let $\Xi^{(n)}$ denote the $n^{th}$ derivative of $\Xi$.  Now $\lambda\mapsto\Xi^{(n)}(z,\nu)(\lambda)$ is analytic on the open disk and extends smoothly to the boundary.  The analyticity guarantees that $$\left|{\partial\over\partial r}\Xi^{(n)}(z,\nu_0)(r)|_{r=r_0}\right |=\left|{\partial\over\partial (i\nu)}\Xi^{(n)}(z,\nu)(r_0)|_{\nu=\nu_0}\right |,$$
where $i$ is the imaginary unit.  For fixed $n$, the right side is uniformly bounded in $r_0$ for $r_0$ near $1$, since $\Xi$ is $C^\infty$, so the left side is also.  Thus $\Xi^{(n)}(z,\nu)(r)\rightarrow\Xi^{(n)}(z,\nu)(1)$ uniformly in $(z,\nu)$ as $r\rightarrow 1$.  Now let $\rho$ be the defining function for $K$ prescribed after (3).  Choose $R<1$ so large that for $r\geq R$, the image $(z,\nu,r)\mapsto(z,\Xi(z,\nu)(r))$ lies in the domain where $\rho$ is strictly hypoconvex and where projection to $K_z$ is nonsingular.  Let $\pi(\Xi(z,\nu)(r))$ be the point on $K_z$ nearest $\Xi(z,\nu)(r)$; this is $C^\infty$ in $(z,\nu,r)$ for $r$ near $1$.  For $\epsilon>0$ and $R\leq r<1+\epsilon$, define $$\Xi_1(z,\nu)(r)={(1-r)(\Xi(z,\nu)(R))+(r-R)\pi(\Xi(z,\nu)(r))\over 1-R}.$$    Then if $R$ is near enough to $1$ and $\epsilon $ close enough to $0$, the mapping $(z,\nu,r)\mapsto(z,\Xi_1(z,\nu)(r))$ is a $C^\infty$ diffeomorphism in a neighborhood of the set where $R\leq r\leq 1$.  For fixed $r$, the image of $(z,\nu)\mapsto(z,\Xi_1(z,\nu)(r))$ satisfies (1); this requires that $\epsilon$ be small enough, that $R$ be close enough to $1$, and makes use of the fact that all derivatives of $\Xi(\cdot,\cdot)(r)$ converge uniformly as $r\rightarrow 1$. Since $(z,\nu,r)\mapsto(z,\Xi_1(z,\nu)(r))$ is a diffeomorphism, we may define $\uk(z,w)$ to be the square of the $r$ coordinate of the inverse.  The functions $\ui,\uk$ coincide on the set where $\ui=R^2$.  Let $\um=f\circ\uk$, where $f$ is an increasing $C^\infty$ diffeomorphism of $[R^2,1]$ to itself such that the derivative of $f$ at $R^2$ is so large that the gradient $|D_w\um|$ on the set where $\ui=R^2$ is greater than $|D_w\ui|$.  If we define $\ul$ to equal $\ui$ where $\ui\leq R^2$ and $\um$ where $\ui\geq R^2$, then $\ul$ would satisfy all the properties we wish $\um$ to satisfy except that it is not smooth on ${\Cal S}:=\{(z,w)\in\Gamma\times\C^n:\ul(z,w)=R^2\}$.  We correct this problem by mollifying $\ul$ near ${\Cal S}$.  We do this by applying Lemma \lemlocmollify, with $\rho_1=\ui-R^2+1$, $\rho_2=\ul-R^2+1$, $K={\Cal S}$, on a small neighborhood $N$ of ${\Cal S}$.  This produces the function $\un$.

Next we modify $\un$ to obtain a function whose level sets near the zero set have spherical fibers over the circle.

Choose $P$ so small that the image $(z,\nu)(r)\mapsto(z,\Xi(z,\nu)(r))$ has strictly convex fibers over the circle for $r\leq P$ (possible by Lemma \lemdelta.)  
Choose $T<P$ so small that the absolute value of the $w$-coordinate of $\Xi(z,\nu)(P)$ is larger than $T$ for all $z,\nu$.  Now define a function $u_7^*(z,w)$ as follows.  The ray from $(z,0)$ through $(z,w)$ pierces the sphere where $|w|^2=T^2$ and the set where $u_6(z,\cdot)$ has the value $P^2$.  Let us define $u_7^*$ to be $T^2$ where the sphere is pierced, $P^2$ where the $P^2$ level set of $\un$ is pierced, and extend affinely to the whole ray.  Then $u_7^*$ is strictly hypoconvex  and its level sets satisfy (1) where $T^2\leq u_7^*\leq P^2$.  Let $I$ be an increasing diffeomorphism of $[T^2,P^2]$ to itself such that $|D_w(I\circ u_7^*)(z,w)|\geq |w|$ when $|w|=T$ and $|D_w(u_6)(z,w)|\geq |D_w(I\circ u_7^*)(z,w)|$ when $u_6(z,w)=P^2$.  Next define $u_7(z,w)$ to be $u_6(z,w)$, if $u_6(z,w)\geq P^2$; $|w|^2$, if $|w|\leq T$; and $I\circ u_7^*(z,w)$ for all other $(z,w)$.  Then $u_7$ is smooth except where $u_7=T^2$ or $P^2$.  We then locally mollify $u_7$ near these sets to obtain the desired defining function $\uo$ - using Lemma \lemlocmollify.  Then $\uo$ is strictly hypoconvex and satisfies the conditions of the lemma. $\square$

{\bf Lemma \lemdeffuncB .}  {\it Let $K$ satisfy (1).  Then there exists a defining function $\ur$ for $K$ satisfying (4).}

We will also define $K^t\equiv\{(z,w)\in\Gamma\times\C^n:\ur(z,w)=t\}$.

{\it Proof.}  We make use of a transformation provided by Lempert [L5, p.{}882, (5.3)].  Given a strictly hypoconvex, bounded, $C^\infty$ bounded domain $D\subset \C^n$ there exists a ``dual'' $D'\subset\C^n $ which is also strictly hypoconvex, bounded, and $C^\infty$ bounded.  For $\epsilon$ small enough that $K^{1+\epsilon}$ is defined, we construct the dual of $\widehat K^{1+\epsilon}_z$ for every $z$ and call its boundary $H_z$; the $H_z$ fit together to form a new compact set $H$ satisfying (1).  We may construct a function $\uo '$ for $H$ as $\uo$ was for $K$ in Lemma 6.  Consider the mapping $$\eqalign{I:\Gamma\times\C^n&\rightarrow\Gamma\times\C^n\cr
(z,w)&\mapsto(z,{{\partial \uo'\over\partial w_j}(z,w)\over \sum_{j=1}^n w_j{\partial \uo'\over\partial w_j}(z,w)}).}$$  Lempert's work shows that this is a $C^\infty$ homeomorphism from $\{(z,w)\in\Gamma\times\C^n:w\in\widehat H_z\setminus \{0\}\}$ to $\{(z,w)\in\Gamma\times\C^n:w\not\in int\,\widehat K^{1+\epsilon}_z\}$ (the ``outside'' of $\widehat K^{1+\epsilon}$.)  We may then transform the function $\uo'$ to a function $\uq$ on $\{(z,w)\in\Gamma\times\C^n:w\not\in int\,\widehat K^{1+\epsilon}_z\}$: $\uq(z,w)=1/\uo'(I^{-1}(z,w))$.  Then the same work shows that the level sets of $\uq$ are compact sets satisfying (1), $\uq$ is identically $1$ on $K^{1+\epsilon}$, $\uq\geq 1$, $\uq$ is strictly hypoconvex when $\uq>1$, and $\uq\rightarrow\infty$ as $|w|\rightarrow\infty$.  Furthermore, when $t$ is large enough, the set where $\uq=t$ has spherical fibers, because the dual of a ball centered at the origin is another ball centered at the origin.   Taking $\uo$ and $(1+\epsilon)\uq$ together produces a real valued function (call it $\us$) defined on all of $\Gamma\times\C^n$ which is strictly hypoconvex except that it is not necessarily smooth when $\us=1+\epsilon$.  However, a mollification near $K^{1+\epsilon}$ - using Lemma \lemlocmollify again - will produce a $\ur$ satisfying (4), except possibly for (4)(e).  But (4)(e) can easily be obtained by making use of the fact that the level sets of $\ur$ have spherical fibers for large $w$, with radius depending only on $z$.  $\square$

$\S 4$ {\bf Proofs of the theorems.}

{\it Proof of Theorem 3.}  Lemma \lemdeffuncB shows that a $K$ satisfying (1) has a defining function satisfying property (3) on page 679 of [Wh].  Theorem 3 then follows entirely from Theorem 2 of [Wh], except for the part about smooth extension to the boundary of $\Delta$.  This follows from the fact that $K$ is a $C^\infty$ manifold: we use results of \v Cirka [\v Ci] and the remark in the last sentence in section 2 of [Wh, p. 694]. $\square$

{\it Proof of Theorem 2.}    From Lemma \lemdeffuncB, we obtain a defining function for $K$ which satisfies condition (3) on p.{}679 of [Wh]; hence we may apply Theorem 3 of [Wh] to conclude that $\phi$ exists and is unique.  The function $\phi$ extends to be in $C^\infty(\Delta)$ by the \v Cirka result just mentioned.  The only remaining fact that needs to be proven is the statement about the polynomial hull of $K$.  To prove this, note that the hull of $K^t$ increases with $t$; furthermore, $\widehat
K=\cap_{t>1}\widehat K^t$.  By
Theorem \fathull, given any point $(z_0,w_0)$ in $\widehat K$ such
that $|z_0|<1$, we have that $(z_0,w_0)$ lies on the graph of some
analytic $\C^n$ valued mapping $f^t$ whose graph lies in $\widehat
K^t$.  Any local uniform limit $f$ of a subsequence of the $f^t$ must
satisfy $f(z)\in \widehat K_z$ for a.e. $z\in\Gamma$ by [HMa,Corollaries 1,2],
so $f=\phi$.  This shows that the point $(z_0,w_0)$ must lie on the
graph of $\phi$, so the graph of $\phi$ is $\widehat
K\cap\Pi^{-1}(\opdisk)$.  $\square$

{\it Proof of Theorem \hullthm.}  Let $\ur$ be as in Lemma \lemdeffuncB, and let $K^t$ be as before.  Recall the definition of $\gamma_\rho$ from (2).  We consider 3 cases:
(a)$\gamma_{\ur} <1$,(b)$\gamma_{\ur} =1$,(c)$\gamma_{\ur} >1$.  In case
(a), by Theorem \optthm there exists $f\in A(\Delta)^n$ such that
$\ur(z,f(z))<1$ for $z\in \Gamma$.  In fact $f$ may have polynomial
coordinates, by approximation.  Thus $K$ has the form of the $K$ in
Theorem \fathull (after a holomorphic change of variable moving the
graph of $f$ to the graph of the zero function) so Theorem \hullthm
will hold in case (a).  In case (b), we have precisely the
circumstances of Theorem \optthm , so Theorem \hullthm holds.  In case
(c), there exists $t>1$ such that $\gamma_{\ur}=t$.  From
Theorem
\optthm , we obtain a unique $\phi$ such that $\ur(z,\phi(z))\leq
\gamma_{\ur}$.  The graph of this $\phi$ contains the only points
which are in $\widehat K^t\cap\Pi^{-1}(\opdisk)$, so these are the
only points which could possibly be in $\widehat
K\cap\Pi^{-1}(\opdisk)$.  We claim the latter set is empty; if not, it
must contain a point of the graph of $\phi$, so it must contain the
whole graph by a theorem of Oka.  But then the points $\phi(z)$ for
$z\in\Gamma$ must be in $\widehat K_z$.  They are not because for
$z\in\Gamma$, $\phi(z)\in K^t_z$, and $K^t_z\cap\widehat
K_z=\emptyset$ since $t>1$. $\square$

Theorem \hullthm shows that as in the case $n=1$ and as is the
case for $K$ with convex fibers over the circle, the polynomial hull
of a $K$ in $\Gamma\times\C^n$ with strictly hypoconvex fibers is
either empty, a single analytic graph or the union of infinitely many analytic
graphs.

We are grateful to Professor L\' aszl\' o Lempert for his interest,
encouragement and suggestions; and to the referee and editor for their assistance.

We refer the reader to [Wh] for a more complete list of references.


\Refs
\widestnumber\key{HMa}

\ref \key AW \by Alexander, H., and J. Wermer \pages 99--109
\paper Polynomial Hulls with Convex Fibers
\yr 1985 \vol 271
\jour Math. Ann.
\endref

\ref \key BL \by Balogh, Zolt\' an M. and Christoph Leuenberger \pages
915--935
\paper Quasiconformal contactomorphisms and polynomial hulls with convex fibers
\yr 1999 \vol 51, no. 5
\jour Canad. J. Math. 
\endref


\ref \key {\v Ce1} \by \v Cerne, Miran \pages 97--105
\paper Smooth families of fibrations and analytic selections of 
polynomial hulls. 
\yr 1995 \vol 52
\jour Bull. Austral. Math. Soc.
\endref

\ref \key {\v Ce2} \by \v Cerne, Miran \pages 27--45
\paper Maximal plurisubharmonic functions and the polynomial hull of a completely circled fibration 
\yr 2002 \vol 40, no. 1
\jour Ark. Mat.
\endref

\ref \key {\v Ci} \by \v Cirka, E.{}M.{}\pages 291--336
\paper Regularity of boundaries of analytic sets
\yr 1983 \vol 45
\jour Math. USSR Sb.
\endref

\ref \key F1 \by Forstneri\v c, Franc \pages 869--889
\paper Polynomial Hulls of Sets Fibered Over the Circle
\yr 1988 \vol 37
\jour Indiana Univ. Math. J.
\endref

\ref \key F2 \by Forstneri\v c, Franc \pages 97--104
\paper Polynomially convex hulls with piecewise smooth boundaries
\yr 1986 \vol 276
\jour Math. Ann.
\endref

\ref \key HMa \by Helton, J. William and Donald E. Marshall \pages 157--184
\paper Frequency domain design and analytic selections
\yr 1990 \vol 39, no. 1
\jour Indiana Univ. Math. J.
\endref

\ref \key HMe \by Helton, J. William and Orlando Merino \pages 285--287
\paper A fibered polynomial hull without an analytic selection
\yr 1994 \vol 41, no. 2
\jour Michigan Math. J.
\endref

\ref \key H\" o \by H\" ormander, Lars 
\book Notions of Convexity
\publ Birkh\" auser \publaddr Boston
\yr 1994
\endref

\ref \key L1 \by Lempert, L\' aszl\' o \pages 257--261
\paper Holomorphic retracts and intrinsic metrics in convex domains
\yr 1982 \vol 8, no. 4
\jour Anal. Math.
\endref

\ref \key L2 \by Lempert, L\' aszl\' o \pages 341--364
\paper Intrinsic distances and holomorphic retracts
\inbook Complex analysis and applications '81 (Varna, 1981)
\publ Bulgar. Acad. Sci. \publaddr Sofia
\yr 1984
\endref

\ref \key L3 \by Lempert, L\' aszl\' o \pages 427--474
\paper La m\' etrique de Kobayashi et la repr\' esentation des 
domaines sur la boule
\yr 1981 \vol 109, no. 4
\jour Bull. Soc. Math. France
\endref

\ref \key L4 \by Lempert, L\' aszl\' o \pages 515--532
\paper Solving the degenerate complex Monge-Amp\` ere equation with 
one concentrated singularity
\yr 1983 \vol 263, no. 4
\jour Math. Ann.
\endref

\ref \key L5 \by Lempert, L\' aszl\' o \pages 869--885
\paper Symmetries and other transformations of the complex
Monge-Amp\` ere equation
\yr 1985 \vol 52, no. 4
\jour Duke Math. J.
\endref

\ref \key L6 \by Lempert, L\' aszl\' o \pages 43--78
\paper Holomorphic invariants, normal forms, and the moduli space of
convex domains
\yr 1988 \vol 128
\jour Ann. Math.
\endref

\ref \key S1 \by S\l odkowski, Zbigniew \pages 367--391
\paper Polynomial Hulls in {\bf C}$^2$ and Quasicircles
\yr 1989 \vol XVI 
\jour Ann. Scuola Norm. Sup. Pisa Cl. Sci. (4)
\endref

\ref \key S2  \by S\l odkowski, Zbigniew \pages 255--260
\paper Polynomial Hulls with Convex Sections and Interpolating Spaces
\yr 1986 \vol 96, No. 2
\jour Proc. Amer. Math. Soc.
\endref 

\ref \key S3 \by S\l odkowski, Zbigniew \pages 156--176
\paper Polynomial hulls with convex fibers and complex geodesics
\yr 1990 \vol 94, no. 1
\jour J. Funct. Anal.
\endref

\ref \key V \by Vityaev, Andrei E. \pages 161--173
\paper Uniqueness of solutions of a $H^\infty$ optimization
problem and complex geometric convexity
\yr 1999 \vol 9, no. 1
\jour J. Geom. Anal.
\endref

\ref \key Wh \by Whittlesey, Marshall A. \pages 677--701
\paper Polynomial hulls and $H^\infty$ control for a hypoconvex constraint.
\yr 2000 \vol 317
\jour Math. Ann.
\endref

\endRefs
\enddocument